\documentclass[12pt,leqno]{amsart}
\usepackage{lmodern}
\usepackage{textcomp}

\usepackage{a4} 

\usepackage{amsmath}
\usepackage{amssymb}
\usepackage{amsthm}
\usepackage{amscd}
\usepackage{mathdots}

\usepackage{enumerate} 
\usepackage{xcolor}

\usepackage{url}
\usepackage[colorlinks=true,linkcolor=blue,citecolor=blue,urlcolor=blue]{hyperref}
\usepackage[capitalize,nameinlink]{cleveref}
\usepackage{thmtools}

\newcommand{\la}{\langle}
\newcommand{\ra}{\rangle}

\newcommand{\ff}{\mathbb F}
\newcommand{\qq}{\mathbb Q}

\newcommand{\cc}{\mathbb C}
\newcommand{\nat}{\mathbb{N}}
\newcommand{\zz}{\mathbb Z}
\newcommand{\rr}{\mathbb R}

\newcommand{\K}{\mathsf{k}}
\newcommand{\mc}[1]{\mathcal{#1}}

\newcommand{\mg}[1]{#1^{\times}}
\newcommand{\sq}[1]{#1^{\times 2}}
\newcommand{\scg}[1]{\mg{#1}/\sq{#1}}

\renewcommand{\setminus}{\smallsetminus}

\DeclareMathOperator{\newsum}{\mathsf{\Sigma}}
\renewcommand{\sum}{\newsum}
\newcommand{\sums}[1]{\sum\!{#1}^2}
\newcommand{\sos}[2]{\sum_{#1}\!{#2}^2}

\newcommand{\ovl}{\overline}

\DeclareMathOperator{\rk}{\mathsf{rk}}

\newcommand{\wt}[1]{\widetilde{{#1}}}
\newcommand{\mfm}{\mathfrak{m}}

\renewcommand{\inf}{\mathsf{inf}}
\renewcommand{\min}{\mathsf{min}}

\renewcommand{\sup}{\mathsf{sup}}

\newcommand{\car}{\mathsf{char}}
\newcommand{\vf}{\varphi}

\renewcommand{\leq}{\leqslant}
\renewcommand{\geq}{\geqslant}

\renewcommand{\bmod}{\,\,\mathsf{mod}\,\,}

\newcommand{\I}{\mathsf{I}}

\newcommand{\ctd}{\mathsf{cd}_2}

\swapnumbers
\numberwithin{equation}{section}
\newtheorem{thm}[equation]{Theorem}
\newtheorem{prop}[equation]{Proposition}
\newtheorem{cor}[equation]{Corollary}
\newtheorem*{cor*}{Corollary}

\newtheorem*{conj*}{Conjecture}
\newtheorem{lem}[equation]{Lemma}

\newtheorem{qu}[equation]{Question}
\newtheorem*{qu*}{Question}

\theoremstyle{definition}

\newtheorem{ex}[equation]{Example}
\newtheorem{ex*}{Example}
\newtheorem{exs}[equation]{Examples}
\newtheorem{rem}[equation]{Remark}

\theoremstyle{plain}

\title{Sums of squares in function fields over a henselian valued field}
\date{15.09.2026}

\author[K.J.~Becher]{Karim Johannes Becher}
\address{University of Antwerp, Department of Mathematics, Middelheim\-laan~1, 2020 Antwerpen, Belgium}
\email{karimjohannes.becher@uantwerpen.be}

\author[N.~Daans]{Nicolas Daans}
\address{Charles University, Faculty of Mathematics and Physics, Department of Algebra, Sokolov\-sk\' a 83, 18600 Praha~8, Czech Republic}
\curraddr{KU Leuven, Department of Mathematics, Celestijnenlaan 200B, 3001 Heverlee, Belgium}
\email{nicolas.daans@kuleuven.be}

\author[D.~Grimm]{David Grimm}
\address{Universidad de Santiago de Chile, Facultad de Ciencias, Avenida Libertador Bernardo O'Higgins nº 3363, Estaci\'on Central, Santiago, Chile}
\email{david.grimm@usach.cl}

\author[G.~Manzano Flores]{Gonzalo Manzano Flores}
\address{Universidad de Valpara\'iso, Facultad de Ciencias, Instituto de Matem\'aticas, Gran Bretaña 1111, Valpara\'iso, Chile}
\email{gonzalo.manzano@uv.cl}

\author[M.~Zaninelli]{Marco Zaninelli}
\address{University of Pennsylvania, Department of Mathematics, 209 South 33rd\linebreak Street,
Philadelphia, PA 19104-6395, United States of America}
\email{mzan@sas.upenn.edu}

\begin{document}


\subjclass[2020]{11E81, 12D15, 12J10}
\keywords{Pythagoras number, function field in one variable, quadratic form, local-global principle, hereditarily pythagorean field, $2$-cohomological dimension}

\begin{abstract}
The Pythagoras number of a field is the smallest natural number $n$ such that every sum of squares in the field is a sum of $n$ squares. 
Given a nondyadic henselian valued field $(K, v)$, we show that the best bound on the Pythagoras number for function fields of curves over $K$ is essentially the same as for function fields of curves over the residue field $Kv$. As a consequence, we obtain the upper bound $5$ for the Pythagoras number of function fields of curves over any field $K$ for which the rational function field $K(X)$ has Pythagoras number $2$.
The proof uses a reduction to the case where the value group is a subgroup of the real numbers, where then a local-global principle for quadratic forms due to V.~Mehmeti is applied.
\end{abstract}

\maketitle

\section{Introduction} 

There are some remarkable early achievements in the study of sums of squares in fields.
On 17 June 1751 Leonhard~Euler presented to the \emph{Berliner Akademie der Wissenschaften}
his proof that every positive rational number is a sum of four squares (see \cite[p.~375]{BS07}).
In 1906, Edmund~Landau \cite{Lan06} used a result due to David~Hilbert from \cite{Hil93} to show that, in the rational function field in two variables $\rr(X,Y)$, every positive definite rational function is a sum of four squares.
In 1921 Carl Ludwig Siegel \cite{Sie21} showed that in any number field, an element which is positive under all embeddings into $\rr$ is a sum of four squares.
These three historical results have in common that, for certain fields, they bound by four the number of square terms needed to represent an arbitrary (totally) positive element as a sum of squares.

Let $n\in\nat$.
It is not difficult to show  that $1+X_1^2+\dots+X_n^2$ is not a sum of $n$ squares in the polynomial ring $\rr[X_1,\dots,X_n]$. That this element is not even a sum of $n$ squares in the fraction field $\rr(X_1,\dots,X_n)$ was only proven in 1964 by John William Scott~Cassels \cite{Cas64}.
This gave the first evidence that there are fields containing elements that are sums of squares but which cannot be expressed as sums of five squares.

Let $K$ denote an arbitrary field. 
In 1927, 
Emil~Artin and Otto~Schreier \cite{AS27} established a link between the set of sums of squares in $K$ and the possibility of endowing $K$ with a field ordering, that is, a total order relation compatible with addition and with multiplication by positive elements.
The Artin-Schreier Theorem states that $K$ admits a field ordering provided that $-1$ cannot be expressed as a sum of squares in $K$.
At the same time, Artin \cite{Art27} showed that, if $K$ has characteristic different from $2$, then an arbitrary nonzero element of $K$ is positive at every ordering of $K$ if and only if it can be expressed as a sum of squares in $K$.
He further applied this to solve Hilbert's 17th Problem: A polynomial 
$f\in \rr[X_1,\dots,X_n]$ is a sum of squares in $\rr(X_1,\dots,X_n)$ if and only if $f(x)\geq 0$ for all $x\in \rr^n$.
Hilbert \cite{Hil93} had shown this in 1893 for $n=2$.

In 1967, Albrecht~Pfister \cite{Pfi67} provided a quantitative counterpart to Artin's statement by showing that any element in $\rr(X_1,\dots,X_n)$ which is a sum of squares can be expressed as a sum of $2^n$ squares.

These results motivated the introduction of the following numerical field invariant:
We denote by $p(K)$ the smallest natural number $m$ such that every 
sum of squares in $K$ can be expressed as a sum of $m$ squares in $K$, provided that such a number $m$ exists, and otherwise we set $p(K)=\infty$. 
The systematic study of this field invariant was initiated in 1970 by Pfister \cite{Pfi71b}, who later in \cite{Pfi76} introduced the notation $p(K)$ and baptised it the \emph{Pythagoras number of $K$}.

 We denote by $K(X)$ the rational function field in one variable $X$ over $K$.
In 1971, Yves~Pourchet \cite{Pou71} proved that $p(K(X))\leq 5$ holds when $K$ is a number field.
Recently, Olivier Benoist \cite{Ben25} extended this result by showing that  $p(F)\leq 5$ holds for any field $F$ of transcendence degree $1$ over $\qq$.

The techniques developed to date to prove upper and lower bounds for the Pythagoras number of a field are quite complementary, and still very limited as far as lower bounds are concerned. 

The field $K$ is called (\emph{formally}) \emph{real} if $-1$ is not a sum of squares in $K$, and \emph{nonreal} otherwise. It follows by \cite[Theorem 2]{Cas64} that $p(K(X))\geq p(K)+1$ whenever $K$ is real.
In \cite{CEP71} it was established that $p(\rr(X_1,X_2))=4$.
Together with Pfister's above-mentioned bound from \cite{Pfi67}, these results yield that $n+2\leq p(\rr(X_1,\dots,X_n))\leq 2^n$ whenever $n\geq 2$. 
The challenge to determine $p(\rr(X_1,\dots,X_n))$ precisely or to decide at least whether its growth is rather linear or exponential in $n$ has resisted any further attempts so far. 

This brings to the fore the general problem of bounding the growth of the Pythagoras number under field extensions.
In particular, the following question from \cite[Question 8, p.~334]{Lam73} is still wide open:

\begin{qu}\label{Q-main}
Can one bound $p(K(X))$ in terms of $p(K)$? In particular, does $p(K)<\infty$ imply that $p(K(X))<\infty$?
\end{qu}

Note that a positive answer to this question would imply that
$p(K(X_1,\dots,X_n))$ can be bounded in terms of $p(K)$ for all $n\in\nat$.
We now mention a few cases where 
upper bounds for $p(K(X_1,\dots,X_n))$  for all $n\in\nat$ are known.

The following general statement involves the notion of \emph{cohomological $2$-dimen\-sion} of a field, denoted by $\ctd$, which we discuss in \Cref{S:cd2} below.

\begin{prop}[\Cref{R:cd2}, \Cref{C:transcendental-torsfree-p-bound}]\label{P:ctd-pyth-bound}
Assume that $\ctd(K(\sqrt{-1}))<\infty$. 
Let $n\in\nat$ and $i=\ctd(K(\sqrt{-1}))$. 
Then $p(F)\leq 2^{n+i}$ for every field extension $F/K$ of transcendence degree $n$. In particular $p(K(X_1,\dots,X_n))\leq 2^{n+i}$.
\end{prop}

For $i=0$, this statement was obtained in \cite{EW87}. 
This includes the case where $K$ is real closed (e.g.~$K=\rr$), which had first been settled by Ernst Witt in \cite{Witt34} for $n=1$ and by Albrecht Pfister in \cite{Pfi71} for arbitrary $n\in\nat$.
While \Cref{P:ctd-pyth-bound} is essentially well-known, its proof for arbitrary $i\in\nat$ relies on several deep facts, including the Milnor Conjecture,  established by Vladimir~Voevodsky in \cite{Voe03}.
For lack of a convenient reference, we will outline an explanation for \Cref{P:ctd-pyth-bound} in \Cref{S:cd2}, based on those deep facts.

Two recent results show that the exponent $n+i$ in the upper bound given by \Cref{P:ctd-pyth-bound} can be lowered in some cases:
\begin{itemize}    
    \item A result due to Uwe~Jannsen \cite[Corollary~0.7]{Jan16} in 2016 settled the conjecture by Jean-Louis~Colliot-Th\'el\`ene and Uwe~Jannsen from 1991 that, if $K$ is a number field, then $p(K(X_1,\dots,X_n))\leq 2^{n+1}$ for any $n\in\nat\setminus\{0,1\}$. Note that $\ctd(K(\sqrt{-1})) = 2$, so \Cref{P:ctd-pyth-bound} would yield the upper bound $2^{n+2}$. (See also \cite[Sect.~5]{Pfi00}.)
    \item A result due to Olivier~Benoist \cite[Theorem 0.2]{Ben20} implies (using \Cref{T:Pfister-BVG} below) that, if $r\in\nat$ and $K$ is the fraction field of an $(r+1)$-dimensional regular henselian 
    excellent local ring
    with real residue field $K_0$ (e.g.~$K=K_0(\!(t_0,\dots,t_r)\!)$) and $s = \ctd(K_0(\sqrt{-1}))<\infty$, then 
    $p(K(X_1,\dots,X_n))\leq 2^{n+r+s}$. In this case $\ctd(K(\sqrt{-1})) = r+s+1$, so \Cref{P:ctd-pyth-bound} would yield the upper bound $2^{n+r+s+1}$.
\end{itemize}

Generally, neither $p(K)$ nor $p(K(X))$ is closely related to $\ctd(K(\sqrt{-1}))$:  
For $r\in\nat$, assuming that $\car(K)\neq 2$, passing from $K$ to the iterated Laurent series field 
$K'=K(\!(t_1)\!)\dots (\!(t_r)\!)$, one has $\ctd(K')=\ctd(K)+r$,
while $p(K'(X))$ and $p(K(X))$ are bounded by the same $2$-power, by \cite[Theorem 4.14]{BGVG14}.

As \Cref{Q-main} remains open, it is worthwhile to consider the following variant.

\begin{qu}
\label{Q-main2}
Can one bound $p(K(X,Y))$ in terms of $p(K(X))$?
\end{qu}

Obviously, a positive answer to \Cref{Q-main} would imply a positive answer to \Cref{Q-main2}.
However, one might expect a closer link between $p(K(X,Y))$ and $p(K(X))$ than between $p(K(X))$ and $p(K)$, and \Cref{Q-main2} might therefore be easier to answer than \Cref{Q-main}.

The problem of bounding $p(K(X,Y))$ is closely related to bounding $p(F)$ for function fields in one variable over $K$.
By a \emph{function field in one variable over $K$} we mean a finitely generated field extension of $K$ of transcendence degree $1$.
Our main result, from which, together with \Cref{T:main-2}, we will be able to derive a positive answer to \Cref{Q-main2} in a significant special case, is the following.

\begin{thm}[\Cref{T:p-tilde-henselian-equality} \& \Cref{R:main-link}]\label{T:main}
Let $v$ be a henselian valuation on $K$ with residue field $Kv$ of characteristic different from $2$.
Let $d\in\nat$, not a power of $2$.
Then $p(E) \leq d$ for every function field in one variable $E/Kv$ if and only if
 $p(F) \leq d$ for every function field in one variable $F/K$.
\end{thm}
The proof of \Cref{T:main} crucially relies on a new local-global principle for isotropy of quadratic forms due to Vler\"e~Mehmeti \cite{Meh19}, which we will explain in \Cref{S:Mehmeti}. 
\Cref{T:Mehmeti} states a new variant of this local-global principle directly applicable to our purposes.
In the setting where $v$ is a complete valuation with value group $\zz$, that local-global principle had been established earlier by Jean-Louis~Colliot-Thélène, Raman Parimala and Suresh Venapally \cite{CTPS12}.

Several of our results presented here are extensions of results that were obtained in \cite{BGVG14}, based on \cite{CTPS12}, in that restricted setting.
In particular, \Cref{T:main} extends \cite[Theorem 6.8]{BGVG14}.
However, the way we apply the result from \cite{Meh19} to obtain \Cref{T:main} does not simply extend the line of argument from \cite{BGVG14}, but several extra steps become necessary to eliminate conditions on the value group.
Our new technique of applying the local-global principle  from \cite{Meh19} (in its form presented in \Cref{T:Mehmeti}), following a reduction to the rank-$1$ case, has since been used to prove other results for function fields in one variable over a henselian valued field: in \cite{Daa24} on linkage of Pfister forms, in \cite{BDM25} on the $u$-invariant, and in \cite{AB26} on a decomposition property for torsion forms.

By a well-known reduction argument (\Cref{T:Pfister-BVG} below),  \Cref{T:main} leads to part $(b)$ of the following statement, while part $(a)$ is more elementary.

\begin{cor}[\Cref{C:rff1var-2powbound}, \Cref{C:hens-rafufi2var}]\label{MT}
    Let $v$ be a henselian valuation on $K$ with residue field $Kv$ of characteristic different from $2$. Let $k\in\nat$.
\begin{enumerate}[$(a)$]
    \item If $p(Kv(X))\leq 2^k$, then $p(Kv(X))\leq p(K(X))\leq 2^k$.
    \item If $p(Kv(X,Y))\leq 2^k$, then $p(Kv(X,Y))\leq p(K(X,Y))\leq 2^k$.
\end{enumerate} 
\end{cor}

The question, motivated by $(a)$, whether $p(K(X))=p(Kv(X))$ holds, remains  open. This corresponds to \cite[Conjecture 4.16]{BGVG14}.
In \Cref{nthv-rafufi-pyth}, we show that
\[p(K(X))=\sup\{p(\ell(X))\mid \ell/Kv\mbox{ finite field extension}\}\]
if the henselian valuation $v$ is non-trivial. 
This extends \cite[Cor.~6.9]{BGVG14}, where this was shown in the case where $v$ has value group $\zz$.

Closely related is the question whether $p(L(X))\leq p(K(X))$ holds for any finite field extension $L/K$. This is \cite[Conjecture 4.15]{BGVG14}.
In \Cref{C:pyth-rafufi-finite-extension} we obtain a positive answer to the latter question when $K$ admits a non-trivial henselian valuation of residue characteristic different from $2$.

\Cref{MT} leads us to speculate about the general relation, for a henselian valuation $v$ on $K$ with residue field $Kv$, between the two 
functions $\nat\to\nat\cup\{\infty\}$ given by $n\mapsto p(K(X_1,\dots,X_n))$ and by $n\mapsto p(Kv(X_1,\dots,X_n))$, respectively.

Real fields $K$ with $p(K(X))=2$ form a well-studied class of fields.
They are known as \emph{hereditarily pythagorean fields}. 
In \Cref{Section HPF}, we revisit some properties of these fields, including in particular a valuation-theoretic characterisation due to
Ludwig~Br\"ocker \cite{Bro76} (see \Cref{T:HP-val}). 
It yields on such a field $K$ a henselian valuation $v$ such that  $p(E) \leq 4$ for every function field in one variable $E/Kv$.
As a consequence of this characterisation, of \Cref{T:main} (applied with $d=5$) and of \Cref{MT}, we will obtain:

\begin{thm}[\Cref{C:final}, \Cref{R:final}]\label{T:main-2}
Assume that $p(K(X))= 2$.
Then $p(F) \leq 5$ for every function field in one variable $F/K$ and $p(K(X, Y)) \leq 8$.
\end{thm}

We do not show the optimality of either of these bounds. 
Instead, we believe that, under the same hypothesis, one should have that $p(K(X,Y))\leq 4$ and $p(F)\leq 3$ for every function field in one variable $F/K$; see \Cref{Q:final1}. 

\section{Preliminaries} 
\label{sect:preliminaries} 

Let $K$ be a commutative ring. 
We denote by $\car(K)$ its characteristic, by $\mg{K}$ the multiplicative group of $K$ and by $\sq{K}$ the subgroup of squares in $\mg{K}$.
For $n\in\nat$, we denote by $\sos{n}{K}$ the set of elements of $K$ which can be written as sums of $n$ squares in $K$,
and we set \[\sums{K}=\bigcup_{n\in\nat}\sos{n}{K}.\]

The \emph{level} and the \emph{Pythagoras number of $K$}, denoted respectively by $s(K)$ and $p(K)$, are given by
\begin{align*}p(K) & \,=\,  \inf \{n\in\nat\mid \sums{K}=\sos{n}{K}\},\\
s(K) & \,=\,  \inf \{n\in\nat\mid -1\in\sos{n}{K}\}\,,
\end{align*}
where the infima are taken in $\nat\cup\{\infty\}$.
Note that, if $\car(K)=2$, then we have $p(K)=s(K)=1$.

In the sequel we assume that $K$ is a field.
Note that $K$ is nonreal if and only if $s(K)<\infty$.
If $K$ is nonreal, then one has $s(K)\leq p(K)\leq s(K)+1$ (see e.g.~\cite[Chap.~XI, Theorem 5.6 (2)]{Lam05}) and $s(K)$ is a power of $2$, by a famous result due to Albrecht~Pfister \cite[Chap.~3, Theorem 3.1]{Pfi95}. In particular, in this case $p(K)$ is of the form $2^n$ or $2^n+1$ for some $n\in\nat$.
On the other hand, Detlev~Hoffmann showed in \cite{Hof99} that any positive integer is the Pythagoras number of some real field.

As in \cite[Sect.~4]{BGVG14}, we define
\[p'(K)=\left\{\begin{array}{ll} p(K) & \mbox{ if $K$ is real or $\car(K)=2$},\\ s(K)+1 &\mbox{ if $K$ is nonreal and $\car(K)\neq 2$}.\end{array}\right.\]
Note that $p(K)\leq p'(K)\leq p(K)+1$.
The advantage of the auxiliary field invariant $p'$ compared to $p$ lies in its better behaviour with respect to henselian valuations.
In particular, we have 
\[p'(K)=p(K(\!(t)\!))=p'(K(\!(t)\!)).\]
This is obvious if $\car(K)=2$ and corresponds to \cite[Cor.~4.5]{BGVG14} if $\car(K)\neq 2$.

For a valuation $v$ on $K$, we refer to the pair $(K,v)$ as a \emph{valued field},
we denote by $\mathcal{O}_v$ the valuation ring of $v$, by $\mfm_v$ its maximal ideal, 
by $Kv$ the residue field $\mathcal{O}_v/\mathfrak{m}_v$ and by $vK$ the value group of $v$.
Our standard reference for the theory of valued fields is \cite{EP05}.

\begin{lem}\label{L:valuation-residue-level}
Let $v$ be a valuation on $K$ and 
 $n\in\nat$ such that $1\leq n\leq s(Kv)$.
Then $v(x_1^2+\dots+x_n^2)=2\min\{v(x_i)\mid 1\leq i\leq n\}$ for any $x_1,\dots,x_n\in K$.
In particular $\mc{O}_v\cap\sos{n}{K}=\sos{n}{\mc{O}}_v$.
\end{lem}
\begin{proof}
Since $n\leq s(Kv)$, no sum of at most $n$ nonzero squares in $Kv$ is equal to zero.
As $Kv=\mc{O}_v/\mfm_v$, we obtain for any $x_1\in\mg{\mc{O}}_v$ and $x_2,\dots,x_n\in \mc{O}_v$ that
$x_1^2+\dots+x_n^2\in\mg{\mc{O}}_v$. This readily yields the statement.
\end{proof}

A valuation $v$ on $K$ is called \emph{henselian} if it extends uniquely to any finite field extension of $K$, or equivalently, if every simple root of a polynomial $f\in\mc{O}_v[X]$ in $\mc{O}_v/\mfm_v$ is the residue modulo $\mfm_v$ of a root of $f$ in $\mc{O}_v$. 
The equivalence of these characterisations corresponds to $(1)\Leftrightarrow (4)$ in \cite[Theorem 4.1.3]{EP05}.

\begin{prop}\label{inequality p'}
Let $v$ be a valuation on $K$.
We have that  $s(K)\geq s(Kv)$ and $p(K)\geq p(Kv)$.  
Furthermore, if $\car(Kv)\neq 2$ and $v$ is henselian, then 
\[s(K)=s(Kv)\quad \mbox{ and }\quad p'(K)= p'(Kv).\]
\end{prop}

\begin{proof}
It follows from \Cref{L:valuation-residue-level} that, for any natural number $n\leq s(Kv)$ and any $x\in\mg{\mc{O}}_v\cap\sos{n}{K}$, one has $x+\mfm_v\in\sos{n}{(Kv)}$.
Using this and the fact that $p(Kv)\leq s(Kv)+1$, one readily sees that $s(K)\geq s(Kv)$ and $p(K)\geq p(Kv)$. 

Assume now that $v$ is henselian and $\car(Kv)\neq 2$.

We consider first the case where $s(Kv)<\infty$.
We set $s=s(Kv)$ and choose $x_1,\dots,x_s\in {\mc O}_v$ such that $x_1^2+\dots+x_s^2 \equiv -1\bmod \mfm_v$. 
Since $v$ is henselian and $\car(Kv)\neq 2$, we can  find $y\in \mc{O}_v$ such that $1+x_1^2+\dots+x_{s-1}^2+y^2 =0$. Hence $s(K)=s=s(Kv)$ and consequently $p'(K)=s+1=p'(Kv)$.

We now consider the case where $s(Kv)=\infty$. Then $s(K)=\infty=s(Kv)$ and in particular $p'(K)=p(K)$ and $p'(Kv)=p(Kv)$. Since $p(K)\geq p(Kv)$, to show equality, we may assume that $p(Kv)<\infty$. Let $r=p(Kv)$.
Consider $x\in \sums{K}\setminus\{0\}$.
Then $v(x)\in 2 vK$, by \Cref{L:valuation-residue-level}, so there exists $t\in \mg K$ such that $t^2x\in\mg{\mc{O}}_v$.
In view of \Cref{L:valuation-residue-level}, it follows that there exist $x_1,\dots,x_{r-1}\in \mc{O}_v$ and $x_r\in\mg{\mc{O}}_v$ such that 
$t^2x\equiv x_1^2+\dots+x_r^2\bmod \mfm_v$.
Since $v$ is henselian and $\car(Kv)\neq 2$, we obtain that 
$t^2x=x_1^2+\dots+x_{r-1}^2+y^2$ for some $y\in \mc{O}_v$. Hence $t^2x\in \sos{r}{K}$ and therefore $x\in\sos{r}{K}$.
This argument shows that $p(K)=p(Kv)$. 

Therefore $p'(K)=p'(Kv)$ holds in every case.
\end{proof}

In the context of the inequalities in \Cref{inequality p'}, we note that it may occur that $p'(K)<p'(Kv)$.

\begin{ex}
Consider the field $R=\bigcup_{n=1}^\infty \rr(\!(t)\!)(\sqrt[n]{t})$.
This field is real closed. 
Hence $p(F)=2$ for every function field in one variable $F/R$, by \cite[Cor.~6.3.4]{Pfi95}.
Consider the function field $F=R(X)(\sqrt{(tX-1)(X^2+1)})$. Note that $F$ is real.
The $t$-adic valuation on $\rr(\!(t)\!)$ extends (in fact, uniquely, for it is henselian) to a valuation on $R$ with value group $\qq$, and
in view of \cite[Cor.~2.2.2 \& Cor.~3.2.3]{EP05}, it extends further (uniquely) to a valuation $v$ on $F$ with $vF=\qq$, $v(X)=0$ and such that the residue of $X$ in $Fv$ is transcendental over $\rr$.
For this valuation, we have $Fv\simeq \rr(X)(\sqrt{-(X^2+1)})$. Since $-1=X^2+(\sqrt{-(X^2+1)})^2$ and $\rr$ is relatively algebraically closed in $\rr(X)(\sqrt{-(X^2+1)})$, we obtain that $s(Fv)=2$. 
Therefore $p'(Fv)=3>2=p(F)=p'(F)$.
\end{ex}

We recall the well-known fact that $p(K(X))$ can be bounded in terms of the Pythagoras numbers of all finite field extensions of $K$.

\begin{thm}[Pfister, Becher, Van Geel] 
\label{T:Pfister-BVG}
For $n\in\nat$, the following are equivalent:
\begin{enumerate}[$(i)$]
\item Either $s(K)\leq 2^n$, or $K$ is real and 
$p(L)< 2^{n+1}$ for every finite real field extension $L$ of $K$. 
\item $s(L)\leq 2^n$ for every nonreal finite field extension $L$ of $K$. 
\item $p(K(X))\leq 2^{n+1}$.
\end{enumerate}
\end{thm}
\begin{proof}
If $K$ is nonreal, then $(i)$--$(iii)$ are trivially equivalent.
See \cite[Theorem 3.5]{BVG09} for the equivalence of these conditions under the assumption that $K$ is real.
\end{proof}

\begin{rem}
For finite field extensions $L/K$ one has 
$p(L)\leq [L:K]\cdot p(K)$; see \cite[Chap.~7, Prop.~1.13]{Pfi95}.
If this bound could be replaced by a uniform bound on $p(L)$ in terms of $p(K)$ for all finite field extensions $L/K$, then a bound on $p(K(X))$ in terms of $p(K)$ would follow by means of \Cref{T:Pfister-BVG},  leading in particular to a positive answer to \Cref{Q-main}.

It is an open question whether there exists a field $K$ with  $p(K(X))>p(K)+2$. 
\end{rem}

The following statement extends \cite[Theorem 4.14]{BGVG14}, where it is formulated for $\zz$-valuations.

\begin{cor}\label{C:rff1var-2powbound}
    Let $v$ be a valuation on $K$.
    Then we have $p(Kv(X))\leq p(K(X))$.
    Furthermore, if $\car(Kv)\neq 2$ and $v$ is henselian, then, for any $n\in\nat$, we have $p(K(X))\leq 2^n$ if and only if $p(Kv(X))\leq 2^n$.
     
\end{cor}
\begin{proof}
By \cite[Cor.~2.2.2]{EP05}, $v$ extends to a valuation $w$ on $K(X)$ with residue field $Kv(X)$.
It follows by \Cref{inequality p'} that $p(Kv(X))= p(K(X)w)\leq p(K(X))$.

Assume now that $\car(Kv)\neq 2$ and $v$ is henselian, and let $n\in\nat$ be such that $p(Kv(X))\leq 2^n$. 
In particular $n\geq 1$, because $1+X^2$ is not a square in $Kv(X)$.
To show that $p(K(X))\leq 2^n$, it suffices by \Cref{T:Pfister-BVG} to check that $s(L)\leq 2^{n-1}$ for every finite nonreal field extension $L/K$.
Consider such an extension $L/K$.
Since $v$ is henselian, it extends uniquely to a valuation on $L$, which we denote by $v'$.
Since $v'$ is henselian, too, we have $s(L)=s(Lv')$, by \Cref{inequality p'}. In particular $Lv'$ is nonreal. 
Since $Lv'/Kv$ is a finite field extension, we obtain by \Cref{T:Pfister-BVG} that $s(Lv')\leq 2^{n-1}$. Hence $s(L)\leq 2^{n-1}$.
\end{proof}

\section{The Pythagoras number of function fields} 
\label{S:cd2}

The study of sums of squares in $K$ is naturally embedded into the study of quadratic forms over $K$.
We use standard notation and terminology from quadratic form theory, for which we refer to \cite{Lam05} and \cite{EKM08}.

In this section, we assume that $\car(K)\neq 2$.
For $n\in\nat$, we denote by $\I^n K$ the $n$th power of the fundamental ideal in the Witt ring of $K$.
If $K$ is nonreal, then the Witt ring of $K$ is an additive torsion group, because the quadratic form $2^r\times \la 1\ra$ is hyperbolic for any $r\in\nat$ with $2^r>s(K)$.
Hence, for any $n\in\nat$, $\I^n K=0$ if and only if $K$ is nonreal and $\I^n K$ is torsion-free.

\begin{prop}\label{P:torsfree-pbound}
    Let $n\in\nat$ be such that $\I^{n+1}K$ is torsion-free. Then $p(K)\leq 2^n$.
\end{prop}
\begin{proof}
Consider $a\in\sums{K}\setminus\{0\}$.
For $r\in\nat$ such that $a\in\sos{2^r}{K}$, the $(r+1)$-fold Pfister form $2^r\times \la 1,-a\ra$ is hyperbolic.
It follows that the $(n+1)$-fold Pfister form $\pi=2^n\times \la 1,-a\ra$ represents a torsion element of $\I^{n+1}K$.
Since $\I^{n+1}K$ is torsion-free, it follows that $\pi$ is hyperbolic.
Hence $a$ is a sum of $2^n$ squares in $K$. This proves that $\sums{K}=\sos{2^n}{K}$.
\end{proof}
\begin{prop}[Elman, Lam, Kr\"uskemper]
\label{ELK} 
Assume that $n\in\nat$ is such that $\I^{n+1}K(\sqrt{-1})=0$. Then $\I^{n+1}K$ is torsion-free and $p(K)\leq 2^n$.
\end{prop}
\begin{proof}
See \cite[Cor.~35.27]{EKM08} for the proof that $\I^{n+1}K$ is torsion-free. 
By \Cref{P:torsfree-pbound}, this implies that $p(K)\leq 2^n$.
\end{proof}

For $n\in\nat$, when $K$ is real, vanishing of $\I^{n+1}K(\sqrt{-1})$ is generally a stronger condition than torsion-freeness of $\I^{n+1}K$. The two conditions are, however, equivalent for fields that are rational function fields over some subfield. This was shown essentially in \cite[Theorem 6.9]{BL13}.

\begin{prop}\label{P:BL-rational-torsfree}
Let $n\in\nat$.
The ideal $\I^{n+1}K(X)$ is torsion-free if and only if $\I^{n+1}K(\sqrt{-1})(X)=0$.
\end{prop}
\begin{proof}
If $\I^{n+1}K(\sqrt{-1})(X)=0$, then $\I^{n+1}K(X)$ is torsion-free, by \Cref{ELK}.
Assume now that $\I^{n+1}K(X)$ is torsion-free.
If $K$ is nonreal, then this means that $\I^{n+1}K(X)=0$, and we obtain by \cite[Cor.~35.13]{EKM08} that $\I^{n+1}K(\sqrt{-1})(X)=0$.
If $K$ is real, then it follows by \cite[Theorem 6.9]{BL13} that $\I^{n+1}K(\sqrt{-1})(X)=0$.
\end{proof}

\begin{rem}\label{R:cd2}
Consider $i\in\nat$ and an arbitrary field $L$ with $\car(L)\neq 2$.
Let $\K_iL$ denote the $i$th Milnor $K$-group modulo $2$ of $L$, as in \cite{EL72k}.
By \cite[Cor.~3.3]{EL72k}, we have that $\I^iL=0$ if and only if $\K_iL=0$.
Let $H^i(L,\mu_2)$ denote the $i$th cohomology group for the trivial action of the absolute Galois group of $L$ on the group $\mu_2=\{\pm 1\}\subseteq \mg L$. 
The so-called Milnor Conjecture as formulated and proven in \cite{Voe03} establishes a natural isomorphism between
$\K_iL$ and $H^i(L,\mu_2)$.
Therefore $H^i(L,\mu_2)=0$  if and only if  $\K_iL=0$,  if and only if $\I^i L=0$.

We denote by $\ctd(K)$ the $2$-cohomological dimension  of $K$ as defined and discussed in \cite[Section 6.5]{NSW08}.
By a \emph{primitive field extension of $K$} we mean a finite field extension generated by a single element and hence isomorphic to $K[X]/(p)$ for an irreducible polynomial $p\in K[X]$. 
The Primitive Element Theorem implies that every finite separable field extension is primitive.
The proof of \cite[Prop.~6.5.11]{NSW08} and the above discussion, respectively, yield the following equalities:
\begin{align*}
    \ctd(K) &\, =\,  \sup\{i\in\nat\mid H^i(L,\mu_2)\neq 0\mbox{ for a primitive field extension $L/K$}\}\\
    & \, = \, \sup\{i\in\nat\mid \I^i L\neq 0\mbox{ for a primitive field extension $L/K$}\}.
\end{align*}

For $i\in\nat$, in view of the exact sequence in \cite[Cor.~21.3]{EKM08}, we have $\I^{i+1}K(X)=0$ if and only if 
$\I^iL=0$ for every primitive finite field extension $L/K$.
Therefore 
\[\ctd(K)  =  \inf\{i\in\nat\mid \I^{i+2}K(X) = 0\}.\]
\end{rem}

\begin{prop}\label{cd2substitute}
Let $m,n\in\nat$.
If $\I^{n+1}K(X)=0$, then $\I^{n+m}F=0$ for every field extension $F/K$ of transcendence degree $m$.
\end{prop}
\begin{proof}
Assume that $\I^{n+1}K(X)=0$. 
We will now use the notion of $2$-cohomological dimension of a field of characteristic different from $2$ and its properties recalled in \Cref{R:cd2}. In these terms, we obtain that $\ctd(K)< n$.
For an arbitrary finitely generated field extension $F/K$ of transcendence degree $m$, it follows by \cite[Theorem 6.5.14]{NSW08} that $\ctd(F)= \ctd(K)+m<n+m$, and consequently that $\I^{n+m}F=0$.
This conclusion obviously carries over to arbitrary field extensions $F/K$ of transcendence degree $m$.
\end{proof}

\begin{cor}\label{C:transcendental-torsfree-p-bound}
Let $n\in\nat^+$ be such that $\I^{n+1}K(\sqrt{-1})(X)=0$.
Let $m\in\nat$ and let $F/K$ be a field extension of transcendence degree $m$. 
Then $\I^{n+m}F$ is torsion-free and $p(F)\leq 2^{n+m-1}$.
\end{cor}
\begin{proof}
By \Cref{cd2substitute}, we have that $\I^{n+m}F(\sqrt{-1})=0$.
It follows by \Cref{ELK} that $\I^{n+m}F$ is torsion-free and $p(F)\leq 2^{n+m-1}$.
\end{proof}

\section{Mehmeti's local-global principle for quadratic forms} 

\label{S:Mehmeti}

In this section, we present with \Cref{T:Mehmeti} a local-global principle for quadratic forms over function fields in one variable over a henselian rank-$1$ valued base field $(K,v)$ with $\car(K)\neq 2$.
Such a statement was established first in \cite[Theorem 3.1]{CTPS12} for the case where $(K,v)$ is complete with value group $vK=\zz$ and $\car(Kv)\neq 2$.
The extension to the setting where $(K,v)$ is only assumed to be a complete rank-$1$ valued field with $\car(K)\neq 2$ was established by V.~Mehmeti in \cite{Meh19}. 
For our main purpose, the study of field invariants, it is crucial to have \Cref{T:Mehmeti} available 
in the henselian setting.
\Cref{T:Mehmeti} will be established by a reduction to the complete case, then invoking \cite[Cor.~3.19 \& Rem.~3.20]{Meh19}.
This reduction relies on valuation-theoretic arguments, which will be detailed. 

We start with some more notation and terminology concerning valuations.
In the sequel, let $(K,v)$ denote a valued field.
For a subfield $L$ of $K$, we denote by $vL$ the value group of $v|_L$, that is, $vL=v(\mg{L})$, and by $Lv$ the residue field of $v|_L$.
We say that $v$ is \emph{discrete} if  $vK\simeq \mathbb{Z}$.
The set of convex subgroups of $vK$ is totally ordered by inclusion. 
We denote by $\rk(v)$ the number of nonzero convex subgroups of $vK$ if there exist only finitely many of them, and otherwise we set $\rk(v)=\infty$.
We refer to $\rk(v)$ as the \emph{rank of $v$}.
(For infinite rank, this is a simplification compared to the definition of the rank of $v$~in \cite[Sect.~2.1]{EP05} as an order type; for the purposes of the present article, no distinction between different infinite ranks would matter.)
We denote by $(K^v,\hat{v})$ the \emph{completion of $(K,v)$} as characterised in \cite[Theorem 2.4.3]{EP05}. 
By \cite[Prop.~2.4.4]{EP05} we have $vK = \hat{v}K^v$ and $Kv = K^v\hat{v}$. When $K^v=K$, we also say that \emph{$(K,v)$ is complete}.
If $v$ is a rank-$1$ valuation, then $v$ induces an ultrametric on $K$, and $K^v$ is given as the completion of $K$ with respect to that ultrametric.
In particular, if  $(K,v)$ is complete and $\rk(v)=1$, then  Hensel's Lemma \cite[Theorem 1.3.1]{EP05} yields that $(K,v)$ is henselian.

We recall that a finitely generated field extension $F/K$ is called \emph{separable} if it is a separable algebraic extension over a purely transcendental extension of $K$, and an arbitrary field extension is called  \emph{separable} if every finitely generated subextension is separable.
A field extension $F/K$ is called \emph{regular} if it is separable and $K$ is relatively algebraically closed in $F$, or equivalently, if for every field extension $L/K$, $F\otimes_KL$ is a domain; see \cite[Chap.~V, {\S}17, N$^\circ$ 5]{BouAlg}.

Given a field extension $F/E$, we say that $E$ is \emph{existentially closed in $F$} if, for every $l,m,n\in \nat$ and $f_1, \ldots, f_l, g_1, \ldots, g_m \in E[X_1, \ldots, X_n]$, the system
\[\begin{cases}
f_i(x_1, \ldots, x_n) = 0 \quad (1\leq i\leq l)\\
g_i(x_1, \ldots, x_n) \neq 0 \quad (1\leq i\leq m)
\end{cases}\]
has a solution in $E^n$ if and only if it has a solution in $F^n$.
This is compatible with the notion of existential closedness from model theory: The above definition describes precisely that, when $F$ and $E$ are considered as structures in the first-order language of rings, $E$ is an existentially closed substructure of $F$; see for example the discussion in \cite[Section 7.1]{Hodges}.

\begin{lem}\label{L:hens-extclo}
Assume that $(K,v)$ is henselian. Let $K'\subseteq K^v$ be a subfield containing $K$ such that $K'/K$ is separable.
Let $F/K$ be a field extension.
Then $F\otimes_K K'$ is a domain and $F$ is existentially closed in the fraction field of $F\otimes_K K'$.
\end{lem}
\begin{proof}
It is shown in \cite[Theorem 5.9]{Kuh14} that under the given hypotheses $K$ is existentially closed in $K'$. In particular, $K$ is algebraically closed in $K'$, and since $K'/K$ is also separable, it follows that  $F\otimes_K K'$ is a domain.
By definition, this says that $F$ and $K'$ are linearly disjoint over $K$ inside the fraction field of $F \otimes_K K'$.
We conclude by \cite[Lemma 7.2]{DDF21} that $F$ is existentially closed in the fraction field of $F\otimes_K K'$.
\end{proof}

The next three lemmas preparing for the proof of \Cref{T:Mehmeti} will only be applied there to valued fields $(K,v)$ of rank $1$.
Formulating them only for the rank-$1$ case would nevertheless not simplify the exposition.
The following lemma is particularly well-known for rank-$1$ valuations (see e.g.~\cite[Sect.~1.2, p.~14]{EP05}).

\begin{lem}\label{L:complete-fin-ext}
Assume that $(K, v)$ is complete. Let $L/K$ be a finite field extension and $w$ an extension of $v$ to $L$.
Then $(L, w)$ is complete.
\end{lem}

\begin{proof}
See \cite[Cor.~6.24]{Kuhlmann}.
The statement can also be obtained by combining \cite[Theorem 4.15, Cor.~24.3 \& Theorem 24.10]{Warner}.
\end{proof}

\begin{lem}\label{L:HenselianToComplete}
Assume that $(K,v)$ is henselian. Let $F/K$ be a function field in one variable.
There exists a complete valued field $(L,w)$ extending $(K,v)$ and a regular function field in one variable $F''/L$ such that $F\subseteq F''$ 
together with an intermediate field $F'$ such that $F$ is existentially closed in $F'$, the extension $F''/F'$ is  purely inseparable, and $[wL : vK]$ and $[Lw:Kv]$ are finite.
\end{lem}

\begin{proof}
Let $K^\ast$ denote the relative algebraic closure of $K$ in $F$.
Then $v$ extends uniquely to a henselian valuation $v^\ast$ on $K^\ast$, and we obtain that $[K^\ast:K]<\infty$.
By \cite[Theorem 3.3.4]{EP05}, we have $[v^\ast K^\ast:vK]\cdot [K^\ast v^\ast: Kv]\leq [K^\ast:K]<\infty$.
We may therefore replace $(K,v)$ by $(K^\ast,v^\ast)$ and assume that $K^\ast=K$.

We now consider the completion $(K^v,\hat{v})$ of $(K,v)$.
Let $K'$ denote a maximal separable extension of $K$ in $K^v$.
Then $K^v/K'$ is a purely inseparable algebraic field extension.
Since $K=K^\ast$ and $K'/K$ is separable, $F\otimes_K K'$ is a domain.
We denote its fraction field by $F'$.
Then $F'/K'$ is a function field in one variable and $K'$ is relatively separably closed in $F'$.
Moreover, $F$ is existentially closed in $F'$, by \Cref{L:hens-extclo}.
Let $F^\ast$ denote the perfect closure of $F'$.
Since $K'\subseteq F'$ and $K^v/K'$ is purely inseparable, $K^v$ may be viewed as a subfield of $F^\ast$, in a unique way. Let $M$ denote the relative algebraic closure of $K'$ in $F^\ast$.
Then $M$ is the perfect closure of $K^v$. 
Consider now the field compositum $F'M$ inside $F^\ast$. 
Clearly $F'M$ is a function field in one variable over $M$. 
Since $M$ is perfect, any field extension of $M$ is separable; see e.g.~\cite[Theorem 26.3]{Mat89}.
Hence, we can write $F'M=M(t,\alpha)$ for some $t,\alpha\in F'M$ such that $t$ is transcendental over $M$ and $\alpha$ is separable  over $M(t)$. 
Then there exists a finite subextension $M_0/K^v$ of $M/K^v$ such that $t,\alpha\in F'M_0$, $\alpha$ is separable over $M_0(t)$ and $M_0(t,\alpha) =F'M_0$. 
Let $F''=M_0(t,\alpha)$ and let $L$ be the relative algebraic closure of $M_0$ in $F''$. 
Then $M_0\subseteq L\subseteq M$, $[L:M_0]<\infty$ and $F''/L$ is separable.
Hence $F''/L$ is a regular function field in one variable. 
Since $F' \subseteq F'' \subseteq F^\ast$, the extension $F''/F'$ is purely inseparable.
Since $K^v\subseteq L\subseteq M$, the extension $L/K^v$ is purely inseparable. 
The valuation $\hat{v}$ on $K^v$ extends to a valuation $w$ on $L$.
Since $[L:K^v]<\infty$, it follows by \Cref{L:complete-fin-ext} that $(L,w)$ is complete.
Moreover, since $vK=\hat{v}K^v$ and $K^v\hat{v}=Kv$, we have that $[wL:vK]\cdot [Lw:Kv]\leq [L:K^v]<\infty$,  by \cite[Theorem 3.3.4]{EP05}.
\end{proof}

A subgroup of an ordered abelian group $(G,\leq)$ is called \emph{cofinal} if it is not contained in any proper convex subgroup of $(G,\leq)$.

\begin{lem}\label{L:HenselianToComplete2}
Let $(K',v')/(K,v)$ be an extension of valued fields such that $(K', v')$ is complete and $vK$ is cofinal in $v'K'$.
Let $F/K$ and $F'/K'$ be function fields in one variable with $F \subseteq F'$ and $FK' = F'$.
Let $w$ be a discrete valuation on $F'$ which is trivial on $F$ and on $K'$.
Then there exists a valuation $w_0$ on $F$ extending $v$ such that $F^{w_0}$ embeds into $F'w$ as an extension of $F$, and $w_0F$ embeds into $v'K'$.
\end{lem}
\begin{proof}
Since $w$ is trivial on $F$, the field $F$ embeds over $K$ into the residue field $F'w$.
Since $F'/K'$ is a function field in one variable, $F'w$ is a finite extension of $K'$. 
Hence $v'$ extends to a valuation $\tilde{v}$ on $F'w$.
It follows by \Cref{L:complete-fin-ext} that 
$(F'w,\tilde{v})$ is complete.
Let $w_0 = \tilde{v}\vert_{F}$ and 
$e=[\tilde{v}(F'w):v'K']$.
Since $e\leq [F'w:K']<\infty$,
multiplication by $e$ yields an embedding of $\tilde{v}(F'w)$ into $v'K'$.
This implies that $w_0F$ embeds into $v'K'$ and further that  $v'K'$ is cofinal in $\tilde{v}(F'w)$.
Since $vK$ is contained in $w_0F$ and assumed cofinal in $v'K'$, we have that $w_0F$ is cofinal in $\tilde{v}(F'w)$.
It follows that the embedding of valued fields $(F, w_0) \to (F'w, \tilde{v})$ is also an embedding of topological spaces.
We conclude by the universal property of the completion \cite[Theorem 2.4.3]{EP05}, that $F^{w_0}$ embeds over $F$ into $F'w$.
\end{proof}

We are now ready to establish our extension of Mehmeti's local-global principle from \cite{Meh19}.

\begin{thm}\label{T:Mehmeti}
Let $(K,v)$ be a henselian rank-$1$ valued field with $\car(K)\neq 2$. 
Let $F/K$ be a function field in one variable and let $\vf$ be an anisotropic quadratic form over $F$ of dimension at least $3$.
Then there exists a rank-$1$ valuation $w$ on $F$ such that either $w|_K=v$ or $w|_K=0$, 
and such that $\vf$ is anisotropic over the completion $F^w$. 
Moreover, if $v$ is discrete and $\car(Kv)\neq 2$, then $w$ can be chosen to be discrete.
\end{thm}
\begin{proof}
Assume first that $(K,v)$ is complete and $F/K$ is regular.
In this case, the main part of the statement follows by \cite[Cor.~3.19 \& Rem.~3.20]{Meh19}.
Assuming in addition that $v$ is discrete and $\car(Kv)\neq 2$, one may conclude by \cite[Theorem 3.1]{CTPS12} for the entire claim, as any discrete valuation on $K$ is equivalent to $v$ (see e.g.~\cite[Prop.~2.2]{BGVG14}).

In the general situation, assuming $(K,v)$ to be henselian, by \Cref{L:HenselianToComplete} there exist a complete valued field $(L, v')$ extending $(K, v)$, a regular function field in one variable $F''/L$ with $F \subseteq F''$ and a subfield $F'$ of $F''$ containing $F$ such that $F$ is existentially closed in $F'$, $F''/F'$ is purely inseparable and $[v'L : vK] < \infty$, whereby $\rk(v'L)=\rk(vK)=1$.
If $v$ is discrete, then so is $v'$.

Since $F$ is existentially closed in $F'$, $\vf$ remains anisotropic over $F'$.
Since $K$ is of characteristic different from $2$ and $F''/F'$ is purely inseparable, $F''/F'$ is a direct limit of finite extensions of odd degree, so it follows by \cite[Chap.~VII, Theorem 2.7]{Lam05} that
$\vf$ is anisotropic over $F''$.
Hence, by the above, there exists a rank-$1$ valuation $w'$ on $F''$ such that $\vf$ is anisotropic over $(F'')^{w'}$, either $w'|_{L}=v'$ or $w'|_L=0$, and further such that $w'$ is discrete in the case where $\car(Kv)\neq 2$ and $v$ is discrete.

If $w'\vert_F \neq 0$, then we set $w = w'\vert_F$ and obtain that $F^w$ embeds into $(F'')^{w'}$, whereby $\varphi$ is anisotropic over $F^w$.
Moreover, if $v$ is discrete and $\car(Kv)\neq 2$, then $w'$ is discrete, and hence the restriction $w'|_F$ is discrete unless it is trivial. 

It remains to consider the case where $w'\vert_F = 0$. Then $w'|_K=0\neq v$, whereby $w'|_L=0$.
As $(K, v)$ and $(L, v')$ are both rank-$1$ valued fields, $vK$ is cofinal in $v'L$.
By \Cref{L:HenselianToComplete2}, there exists a (necessarily rank-$1$) valuation $w$ on $F$ extending $v$ such that $F^w$ embeds into $(F'')w'$ over $F$ and such that, furthermore, $wF$ embeds into $v'L$, and hence into $vK$.
Recall that $w'$ is a rank-$1$ valuation.
Since $\varphi$ is defined over $F$ and anisotropic over $(F'')^{w'}$, it follows by  \cite[Theorem 7]{Dur48} that $\varphi$ is anisotropic over $F''w'$. 
Since $F^w$ embeds into $(F'')w'$ over $F$, we conclude that $\varphi$ is anisotropic over $F^w$.
Since $wF$ embeds into $vK$, we have further that, if $v$ is discrete, then so is $w$.
\end{proof}

\begin{rem}
Based on \cite[Cor.~3.18 \& Rem.~3.20]{Meh19},  a variant of \Cref{T:Mehmeti} covering the case where $\car(K)=2$ is obtained in \cite[Theorem 7.1]{BDD23}, assuming that $(K, v)$ is complete and $\varphi$ is a non-degenerate quadratic form. We do not expect this version to extend to the case where $(K,v)$ is henselian and $\car(K)=2$ without an extra hypothesis.
\end{rem}

\section{The Pythagoras number of function fields in one variable} 

In view of our goal to give upper bounds on the Pythagoras numbers of function fields in one variable over $K$, we define 
\[\wt{p}(K)=\sup\{p'(F)\mid F/K\mbox{ function field in one variable}\}\,.\]

\begin{ex}\label{E:Pop}
Benoist \cite{Ben25} proved that $\wt{p}(K)\leq 5$ holds for any number field $K$.
Note further that $\wt{p}(\qq)=p(\qq(X))=5$.
\end{ex}

\begin{prop}\quad
\label{P:p-tilde}
\begin{enumerate}[$(a)$]
    \item\label{p-tilde-nr} If $K$ is nonreal, then $\wt{p}(K)=p'(K)=p(K(X))$.
    \item\label{p-tilde-s>1} If $s(K)>1$, then $\wt{p}(K)\geq 3$ and $\wt{p}(K)$ is not a power of $2$.
    \item\label{p-tilde-rafufi2} For $n\in\nat ^+$, we have $\wt{p}(K)<2^{n+1}$ if and only if $p(K(X,Y))\leq 2^{n+1}$.
\end{enumerate}
\end{prop}
\begin{proof}
$(\ref{p-tilde-nr})$ Assume that $K$ is nonreal.
If $\car(K)=2$, then $\wt{p}(K)=p'(K)=p(K(X))=1$.
Assume now that $\car(K)\neq 2$.
For any field extension $F/K$, we have $s(F)\leq s(K)<\infty$ and therefore $p(F)\leq p'(F) \leq s(F)+1\leq s(K)+1\leq p'(K)$.
This readily yields the inequalities $p(K(X))\leq \wt{p}(K)\leq p'(K)$.
On the other hand, $p(K(X))=s(K)+1=p'(K)$, by \cite[Chap.~VII, Prop.~1.5]{Pfi95}. 
Together, this shows the statement.

$(\ref{p-tilde-s>1})$\, 
Assume that $s(K)>1$.
Then $s(K(X)(\sqrt{-(1+X^2)}))=2$ and therefore $\wt{p}(K)\geq p'(K(X)(\sqrt{-(1+X^2)}))=3$.
It remains to show, for any integer $r\geq 2$, that if $\wt{p}(K)\geq 2^r$, then $\wt{p}(K)>2^r$.
So suppose that there exists a function field in one variable $F/K$ with $p'(F)\geq 2^r$.
Note that $2^r-1$ is not a $2$-power, as $r\geq 2$. 
In particular $2^r-1\neq s(F)$. Hence $p(F)\geq 2^r$. We fix an element $x\in \sos{2^r}{F}\setminus \sos{2^r-1}{F}$ and set $F'=F(\sqrt{-x})$. 
Then $s(F')=2^r$ (see \cite[Chap.~XI, Theorem 2.7]{Lam73} or \cite[Chap.~XI, p.~421, Ex.~5]{Lam05}), whereby $\wt{p}(K)\geq p'(F')=2^r+1$.

$(\ref{p-tilde-rafufi2})$\,  This follows by \Cref{T:Pfister-BVG}.
\end{proof}

\begin{cor}\label{C:2powerheadache}
Let \[d=\sup\{p(F)\mid F/K\mbox{ function field in one variable}\}\in\nat\cup\{\infty\}.\]
If $d$ is a $2$-power and $-1\notin\sq{K}$, then $\wt{p}(K)=d+1$, and otherwise $\wt{p}(K)=d$.
\end{cor}
\begin{proof}
If $\car(K)=2$, then $d=1=\wt{p}(K)$.
Assume now that $\car(K)\neq 2$.
If $-1\in\sq{K}$,
then $p(F)=2=p'(F)$ for every function field in one variable $F/K$, whereby 
$d=2=\wt{p}(K)$. We assume now that $-1\notin\sq{K}$.

It follows from the definitions that $d\leq \wt{p}(K)\leq d+1$.
If $d$ is a $2$-power, then $\wt{p}(K)\neq d$, by \Cref{P:p-tilde}, and hence $\wt{p}(K)=d+1$.
Assume now that $d$ is not a $2$-power.
In particular $s(F)\neq d$ for any nonreal field $F$.
Hence, for any function field in one variable $F/K$, we have 
$p'(F)=p(F)\leq d$ or $p'(F)=s(F)+1\leq d$.
Therefore $\wt{p}(K)\leq d$, whereby $\wt{p}(K)=d$.
\end{proof}

\begin{cor}\label{C:fufi2-p4}
We have  $p(K(X,Y))=4$ if and only if $K$ is real and $\wt{p}(K)=3$.
\end{cor}
\begin{proof}
By \Cref{P:p-tilde}~$(\ref{p-tilde-rafufi2})$, we have  $p(K(X,Y))\leq 4$ if and only if $\wt{p}(K)\leq 3$.
If $\car(K)=2$, then $K$ is nonreal and $p(K(X,Y))=1\neq 4$.
If $K$ is nonreal and $\car(K)\neq 2$, then $p(K(X,Y))=s(K)+1\neq 4$.
Assume now that $K$ is real.
Then $\wt{p}(K)\geq 3$ by \Cref{P:p-tilde}~$(\ref{p-tilde-s>1})$.
Furthermore, the polynomial 
\[(1+X^2-2X^2Y^2)^2 +(XY^2-X^3Y^2)^2+(XY-X^3Y)^2+ (X^2Y-X^4Y)^2\]
is not a sum of three squares in $K(X,Y)$; see~\cite[Section~3]{Gri15}. 
Therefore we have $p(K(X,Y))\geq 4$. 
\end{proof}

\begin{lem}\label{L:alg-p-tilde}
For every algebraic field extension $L/K$, one has $\wt{p}(K)\geq p'(L)$.
\end{lem}
\begin{proof}
It suffices to show this for the case where $L/K$ is a finite field extension.
In that case, $L(X)/K$ is a function field in one variable, whereby $\wt{p}(K)\geq p'(L(X))$.
On the other hand, we obtain by \Cref{inequality p'} that $p'(L(X))\geq p'(L)$, because $L$ is the residue field of a 
 valuation on $L(X)$.
\end{proof}

\begin{lem}\label{L:residue-fufi-lift}
Let $v$ be a valuation on $K$ with $\car(Kv)\neq 2$. Let $F_0/Kv$ be a function field in one variable.
Then there exists a function field in one variable $F/K$ and an extension $w$ of $v$ to $F$ such that $Fw$ is $Kv$-isomorphic to $F_0$ and $s(F)=s(F_0)$.
\end{lem}
\begin{proof}
We fix an element $\xi\in F_0$ which is transcendental over $Kv$.
Let $E=K(X)$. By \cite[Cor.~2.2.2]{EP05}, $v$ extends uniquely to a valuation $\wt{v}$ on $E$ with $\wt{v}(X)=0$ and such that the residue $\xi$  of $X$ in $E\wt{v}$ is transcendental over $Kv$, and we obtain that 
$\wt{v}E=vK$ and $E\wt{v}=Kv(\xi)\subseteq F_0$.
Since $F_0/Kv(\xi)$ is a finite field extension, it follows by \cite[Prop.~2.5]{BGVG14} that there exists a finite field extension $F'/E$ with $[F':E]=[F_0:E\wt{v}]$ together with an extension $w'$ of $\wt{v}$ to $F'$ whose residue field $F'w'$ is equal to $F_0$.

It follows that $s(F')\geq s(F_0)$.
If $s(F')=s(F_0)$, then we take $F=F'$ and $w=w'$.
Assume now that $s(F')>s(F_0)$.
In particular, $F_0$ is nonreal.
Let $s=s(F_0)$.
If $s=1$, then we set $f=1$, otherwise we choose $f\in \sos{s}{\mc{O}_{w'}}$
such that $\ovl{f}=-1$ in $F'w'=F_0$.
We set $F=F'(\sqrt{-f})$ and fix an extension $w$ of $w'$ to $F$.
The choice of $f$ implies that 
$Fw=F'w'=F_0$ and $s(F)=s=s(F_0)$, in view of \cite[Chap.~VII, Theorem 3.1]{Lam05}.
\end{proof}

\begin{prop}\label{P:wt-p-geq}
Let $v$ be a valuation on $K$. Then $\wt{p}(K)\geq \wt{p}(Kv)$.
\end{prop}
\begin{proof}
If $\car(Kv)=2$, then $\wt{p}(Kv)=1\leq \wt{p}(K)$.
Assume now that $\car(Kv)\neq 2$.
Consider a function field in one variable $F_0/Kv$.
By \Cref{L:residue-fufi-lift}, there exists a function field in one variable $F/K$ and an extension $w$ of $v$ to $F$ such that $Fw$ is $Kv$-isomorphic to $F_0$ and $s(F)=s(F_0)$.
Using \Cref{inequality p'}, it follows that $p'(F_0)=p'(Fw)\leq p'(F)\leq \wt{p}(K)$.
Since this holds for an arbitrary function field in one variable $F_0/Kv$, the statement follows.
\end{proof}

\begin{thm}\label{T:p-tilde-henselian-equality}
Let $v$ be a henselian valuation on $K$ such that $\car(Kv)\neq 2$.
Then \[\wt{p}(K)=\wt{p}(Kv).\]
\end{thm}

\begin{proof}
By \Cref{P:wt-p-geq}, we have that $\wt{p}(K)\geq\wt{p}(Kv)$.
So it remains to show the other inequality, and to this end, we may assume that $\wt{p}(Kv)<\infty$. 
Note that $s(K)=s(Kv)$.
If $K$ is nonreal, then $\wt{p}(K)=s(K)+1=s(Kv)+1=\wt{p}(Kv)$.
So we may further assume that $K$ is real. Then $Kv$ is real, 
and we obtain that $\wt{p}(K)\geq \wt{p}(Kv)\geq 3$.
Hence, if $\wt{p}(K)=3$, then we have the stated equality.
We may therefore now assume that $\wt{p}(K)>3$.

Note that, if $v$ is the trivial valuation, then $K=Kv$ and thus $\wt{p}(K)=\wt{p}(Kv)$.

Now we consider the case where the rank of $v$ is $1$.
Consider any integer $d\geq 3$ and any function field in one variable $F/K$ such that $p'(F)>d$.
Then either $\sos{}{F}\neq \sos{d}{F}$ or $d\leq s(F)<\infty$. 
If $\sums{F}\neq \sos{d}{F}$ we fix $a\in\sums{F}\setminus\sos{d}{F}$ and set $\vf=d\times\la 1\ra\perp\la -a\ra$, and if $d\leq s(F)<\infty$ we set $\vf=d\times \la 1\ra$.
In either case $\vf$ is an anisotropic quadratic form of dimension at least $3$ over $F$.
Hence, by 
\Cref{T:Mehmeti}, there exists a rank-$1$ valuation $w$ on $F$ whose restriction to $K$ is trivial or equal to $v$ and such that $\vf$ is anisotropic over the completion $F^w$.
By the choice of $\vf$, we obtain that 
$p(F^w)> d$ or $d\leq s(F^w)<\infty$.
Therefore $p'(F^w)> d$.
Since $Fw$ is the residue field of the natural extension of $w$ to $F^w$, we obtain by \Cref{inequality p'} that 
\[p'(Fw)=p'(F^w)> d.\]

If $w$ is trivial on $K$, then $Fw$ is a finite extension 
of $K$, and hence $v$ extends to a henselian valuation $v'$ on $Fw$, whose residue field $(Fw)v'$ is a finite field extension of $Kv$, and in view of \Cref{L:alg-p-tilde}, we obtain that 
\[ \wt{p}(Kv)\geq p'((Fw)v')=p'(Fw)>d.\]
If $w$ extends $v$, then $Fw/Kv$ is either an algebraic field extension or a function field in one variable, and in view of \Cref{L:alg-p-tilde}, we obtain in either case that \[\wt{p}(Kv)\geq p'(Fw)>d.\]

This argument shows that $p'(F)\leq \wt{p}(Kv)$ for all function fields in one variable $F/K$.
We thus have established the equality $\wt{p}(K)=\wt{p}(Kv)$ in the case where $v$ has rank $1$.

In the next step, we extend the statement to the case where $v$ has finite rank.
We denote the rank of $v$ by $n$ and assume that $0<n<\infty$.
We will show the statement by induction on $n$.
Since $n<\infty$, it follows that the value group $vK$ has a unique maximal proper convex subgroup $\Delta$.
Let $v_1$ denote the valuation on $K$ with value group $vK/\Delta$ given by $x\mapsto v(x)+\Delta$.
Then $v_1$ is a rank-$1$ valuation on $K$.
Furthermore $v$ induces a valuation $\wt{w}$ on the residue field $Kv_1$ with value group $\wt{w}Kv_1=\Delta$ and hence of rank $n-1$, and with residue field $(Kv_1)\wt{w}=Kv$.
By \cite[Cor.~4.1.4]{EP05} the valuation $v_1$ on $K$ and the valuation $\wt{w}$ on $Kv_1$ are both henselian.
As $v_1$ has rank $1$, we already know that $\wt{p}(K)=\wt{p}(Kv_1)$.
Furthermore, the induction hypothesis applied to the valuation $\wt{w}$ on $Kv_1$ yields that  $\wt{p}(Kv_1)=\wt{p}(Kv)$, because $Kv=(Kv_1)\wt{w}$. Hence $\wt{p}(K)=\wt{p}(Kv)$.
This proves the statement in the case where $v$ has finite rank.

We now consider the general case, where no assumption is made on the rank of $v$.
We consider a function field in one variable $F/K$. 
If $F$ is nonreal, let $f=-1$, and otherwise, consider an arbitrary element $f\in\sums{F}\setminus\{0\}$.
Let $m\in\nat$ be minimal such that $f \in \sos{m}{F}$.
Recall that $K$ is real, so in particular the smallest subfield of $K$ is $\qq$.
Let $f_1,\dots,f_m\in F$ be such that $f=f_1^2+\dots+f_m^2$.
Let $F_0/\qq(f_1,\dots,f_m)$ be a finitely generated extension such that $F=KF_0$ and $F_0/(K\cap F_0)$ has transcendence degree $1$.
Then $F_0/\qq$ is finitely generated and $f\in\sos{m}{F_0}$.
It follows that $F_0/(F_0\cap K)$ is a function field in one variable and the field extension $(F_0\cap K)/\qq$ is finitely generated.
By \cite[Prop.~3.4.1 and Theorem~3.4.3]{EP05}, it follows that the valuation $v|_{(F_0\cap K)}$ has finite rank.
Let $\Gamma=\{\gamma\in vK\mid \exists n\in\nat^+:n\gamma\in v(F_0\cap K)\}$.
Note that $\Gamma$ is an ordered subgroup of $vK$ containing $v(F_0\cap K)$ and with $\rk(\Gamma)=\rk(v(F_0\cap K))<\infty$.
By Zorn's Lemma, there exists a subfield $K_0$ of $K$ containing $F_0\cap K$ which is maximal for the property that $vK_0\subseteq \Gamma$. 
Then $v(F_0\cap K)\subseteq vK_0\subseteq\Gamma$, and it follows that $\rk(vK_0)=\rk(v(F_0\cap K))<\infty$.
By the choice of $\Gamma$ and $K_0$, it follows that $K_0$ is relatively algebraically closed in $K$.

We now exploit the hypothesis that $(K,v)$ is henselian.
Since $K$ is real, it follows that $Kv$ is real and in particular of characteristic zero.
Since $K_0$ is relatively algebraically closed in $K$, we conclude further that $(K_0,v|_{K_0})$ is henselian, 
by \cite[Cor.~4.1.5]{EP05},
and that $K_0v$ is relatively algebraically closed in $Kv$.
Consider $x\in\mc{O}_v$ and suppose that $x+\mfm_v\in Kv\setminus K_0v$.
Then $x+\mfm_v$ is transcendental over $K_0v$, whereby $vK_0(x)=vK_0$ by \cite[Cor.~2.2.2]{EP05}, which contradicts the choice of $K_0$.
This argument shows that $K_0v=Kv$.
Since $(K_0,v|_{K_0})$ is henselian and $\rk(vK_0)<\infty$, we conclude by the previous case that 
\[\wt{p}(K_0)= \wt{p}(K_0v)= \wt{p}(Kv).\]

Recall that $f\in\sos{m}{F_0}\setminus \sos{m-1}{F}$.
If $F$ is nonreal, then $f=-1$, and we obtain that $m=s(F)=s(F_0K_0)<p'(F_0K_0)\leq \wt{p}(K_0)=\wt{p}(Kv)$, whereby $p'(F)=m+1\leq \wt{p}(Kv)$.
Assume now that $F$ is real. Then so is $F_0K_0$, and hence $m\leq p(F_0K_0)\leq \wt{p}(K_0)=\wt{p}(Kv)$.
Applying this argument for arbitrary $f\in\sums{F}$ and $m\in\nat$ minimal with $f\in\sos{m}{F}$,
we conclude that $p'(F)=p(F)\leq \wt{p}(Kv)$.

Hence we have $p'(F)\leq \wt{p}(Kv)$ for all function fields in one variable $F/K$.
This shows that $\wt{p}(K)\leq \wt{p}(Kv)$.
\end{proof}

\begin{rem}\label{R:main-link}
For $d\in\nat\setminus\{2^i\mid i\in\nat\}$, \Cref{C:2powerheadache} yields that $d\geq \wt{p}(K)$ if and only if $d\geq p(F)$ for every function field in one variable $F/K$. 
Therefore \Cref{T:p-tilde-henselian-equality} yields \Cref{T:main}.
\end{rem}

The following two examples show that \Cref{T:p-tilde-henselian-equality} becomes false if either of the two hypotheses on $v$ is removed.

\begin{ex}
Let $L=\rr(\!(t)\!)$ and let $v$ be the $t$-adic valuation on $L$, which is discrete.
For any subfield $K$ of $L$ properly containing $\rr$, we have $Kv=\rr\not\simeq K$, whereby $v|_K\neq 0$ and thus $vK\simeq \zz$.
Since $L/\rr$ has infinite transcendence degree, we may choose a family $(x_i)_{i\in\nat}$ in $L$ which is algebraically independent over $\rr$.
Let $n\in\nat$ and set $K_n=\rr(x_0,\dots,x_n)$.
Then $\wt{p}(K_n)\geq p(K_n)+1>n$.
On the other hand $\wt{p}(K_nv)=\wt{p}(\mathbb{R})=3$.
\end{ex}

\begin{ex}    
We have $\wt{p}(\qq_2)=p(\qq_2(X))=5$, whereas $\wt{p}(\ff_2)=p(\ff_2(X))=1$. 
More generally, if $K$ carries a henselian valuation $v$ with $\car(Kv)= 2$, then $-7$ is a square in $K$, and it follows that $p'(F)\leq s(F)+1\leq 5$ for any field extension $F/K$, so one can at least conclude that $\wt{p}(K)\leq 5$.
\end{ex}

\begin{cor}\label{C:hens-rafufi2var}
Let $v$ be a henselian valuation on $K$ such that $\car(Kv)\neq 2$. Let $n\in\nat$.
Then $p(K(X,Y))\leq 2^n$ if and only if $p(Kv(X,Y))\leq 2^n$.
\end{cor}
\begin{proof}
The statement is vacuous for $n=0$. 
It follows from \Cref{T:Pfister-BVG} that, for any field $k$, one has $p(k(X,Y))\leq 2$ if and only if $s(k)=1$.
Since $s(Kv)=s(K)$, by \Cref{inequality p'}, this shows the statement for $n=1$.
For $n\geq 2$, the statement follows from \Cref{T:p-tilde-henselian-equality} and \Cref{P:p-tilde}~$(\ref{p-tilde-rafufi2}$).
\end{proof}

\begin{ex}
Let $K$ be a field carrying a henselian valuation $v$ such that the residue field $Kv$ is a number field.
Then $\wt{p}(K)=\wt{p}(Kv)\leq 5$, in view of \Cref{T:p-tilde-henselian-equality} and \Cref{E:Pop}. That is, we have $p(F)\leq 5$ for every function field in one variable $F/K$. By \Cref{P:p-tilde}, we  conclude that $p(K(X,Y))\leq 8$.
These arguments apply for example to the situation where $K$ is the field of Puiseux series over $\qq$ (see \Cref{EX:hp}).
\end{ex}

The following statement extends \cite[Cor.~6.9]{BGVG14}, where it is established under the more restrictive hypothesis that $v$ is discrete and complete.

\begin{thm}\label{nthv-rafufi-pyth}
Let $v$ be a non-trivial henselian valuation on $K$ such that\linebreak $\car(Kv) \neq 2$. 
Then \[p(K(X))=\sup \{p(\ell (X))\mid \ell/Kv\mbox{ finite field extension}\}.\]
\end{thm}

\begin{proof}
Set $m=\sup \{p(\ell (X))\mid \ell/Kv\mbox{ finite field extension}\}$. 
Since $\car(Kv)\neq 2$, we have that $m\geq 2$.
We may assume that $Kv$ is real. 
We will use multiple times in this proof the fact that $p(L) \leq p'(L) \leq p'(L(X))=p(L(X))$ for any field $L$.
    
Let $\ell/Kv$ be a finite field extension.
In view of \cite[Theorem 3.5]{BGVG14}, $v$ extends to a valuation $w$ on $K(X)$ with $K(X)w\simeq \ell(X)$, which yields that $p(K(X))\geq p(\ell(X))$ by \Cref{inequality p'}.
Hence $p(K(X))\geq m$.

It remains to show that $p(K(X))\leq m$. We may assume that $m<\infty$.
The same reduction steps which are used in \Cref{T:p-tilde-henselian-equality} apply here to reduce the problem to the situation where $\rk\, v=1$.
Assume that we are in this case.
    
Consider an arbitrary rank-$1$ valuation $w$ on $K(X)$ such that $\mc{O}_v \subseteq \mc{O}_w$.
We claim that $p(K(X)^w) \leq m$.
    
Assume first that $K\subseteq \mc{O}_w$.
Then $K(X)w/K$ is a finite field extension, and $v$ extends uniquely to a valuation $v'$ on $K(X)w$, which is henselian as well.
Furthermore $(K(X)w)v'\simeq \ell$ for some finite field extension $\ell/Kv$.  
In view of \Cref{inequality p'} and \Cref{P:p-tilde}~$(a)$, it follows that \[p(K(X)^w)\leq p'(K(X)^w)= p'(K(X)w)=p'((K(X)w)v')=p'(\ell)\leq p(\ell(X))\leq m.\]
     
Assume now that $K\not\subseteq \mc{O}_w$. Then $\mc{O}_w\cap K=\mc{O}_v$ and $K(X)w$ is an extension of $Kv$. 
If $K(X)w/Kv$ is algebraic, then $p'(K(X)^w)=p'(K(X)w)\leq m$.
Otherwise, Ohm's Theorem \cite{Ohm83} yields that $K(X)w\simeq\ell(X)$ for some finite field extension $\ell/Kv$.
Consequently,
\[p(K(X)^w)\leq p'(K(X)^w)=p'(K(X)w)=p'(\ell(X))=p(\ell(X))\leq m.\] 

We are now ready to show that $p(K(X))\leq m$.
Consider $a \in \sos{}{K(X)}$ and the quadratic form $\varphi=\sum_{i=1}^m X_i^2-aX_0^2$.
For an arbitrary rank-$1$ valuation $w$ on $K(X)$ with $\mc{O}_v \subseteq \mc{O}_w$, since $p(K(X)^w) \leq m$, we have that $\varphi$ is isotropic over $K(X)^w$.
It follows by \Cref{T:Mehmeti} that $\varphi$ is isotropic over $K(X)$, whence $a \in \sos{m}{K(X)}$.
This argument shows that $\sos{}{K(X)}=\sos{m}{K(X)}$, whereby $p(K(X))\leq m$.
\end{proof} 

\begin{cor}\label{pyth-rafufi-henselian}
Let $(K,v)$ and $(K',v')$ be two non-trivially henselian valued fields with $Kv= K'v'$ and $\car(Kv)\neq 2$. Then $p(K(X))=p(K'(X))$.
\end{cor}
\begin{proof}
This is obvious from \Cref{nthv-rafufi-pyth}.
\end{proof}

As already mentioned in the introduction, it is an open question whether the inequality $p(L(X))\leq p(K(X))$ holds for all finite field extensions $L/K$.
We now obtain a partial positive answer.

\begin{cor}\label{C:pyth-rafufi-finite-extension}
Assume that $K$ carries a non-trivial henselian valuation $v$ with $\car(Kv)\neq 2$.
Then $p(L(X))\leq p(K(X))$ holds for every finite field extension $L/K$.
\end{cor}
\begin{proof}
Set $K'=K(\!(t)\!)$ and let $w$ denote the $t$-adic valuation on $K'$.
Since $w$ is a non-trivial henselian valuation on $K'$ with $K'w=K$, we obtain by \Cref{nthv-rafufi-pyth} that $p(L(X))\leq p(K'(X))$ holds for every finite field extension $L/K$.
    
Let $v'$ denote a valuation on $K'$ obtained as the composition of $w$ with $v$.
Then $v'$ is a non-trivial henselian valuation on $K'$ with $K'v'=Kv$.
Therefore \Cref{pyth-rafufi-henselian} yields that $p(K'(X))=p(K(X))$.
Hence we conclude that $p(L(X))\leq p(K(X))$ for every finite field extension $L/K$.
\end{proof}

\section{Hereditarily pythagorean fields} 

\label{Section HPF}

In this final section we will draw some conclusions from the results of the preceding sections for function fields over hereditarily pythagorean fields.

The field $K$ is called \emph{pythagorean} if $p(K)=1$, and it is called \emph{hereditarily pythagorean} if $K$ is real and $p(L)=1$ for every finite real field extension $L/K$.

As a special case of \Cref{T:Pfister-BVG}, one retrieves a characterisation of hereditarily pythagorean fields by Eberhard Becker in \cite[Chap.~III, Theorem 4]{Bec78}.

\begin{thm}[Becker]
\label{T:Becker}
We have $p(K(X))\leq 2$ if and only if either $-1\in\sq{K}$ or $K$ is hereditarily pythagorean.
\end{thm}
\begin{proof}
This follows from the equivalence $(i)\Leftrightarrow(iii)$ in 
\Cref{T:Pfister-BVG}  for $n=0$.
\end{proof}

We now revisit some crucial properties of hereditarily pythagorean fields concerning their finite extensions, which we will need.

\begin{lem}\label{L:groups}
Let $G$ be a group and $H$ a subgroup of index $2$ of $G$ such that $\sigma^2=1_G$ holds for every $\sigma\in G\setminus H$.
Then $H$ is abelian and every subgroup of $H$ is normal in $G$.
\end{lem}
\begin{proof}
For any $\tau\in H$ and $\sigma\in G\setminus H$, we have $(\sigma\tau)^2=1_G=\sigma^2$, and therefore $\sigma\tau\sigma^{-1}=\sigma\tau\sigma=\tau^{-1}$.
Now fix some $\sigma\in G\setminus H$. For any $\tau,\tau'\in H$ and $\sigma\in G\setminus H$, we obtain that
$\tau'\tau=(\sigma\tau'^{-1}\sigma)(\sigma\tau^{-1}\sigma)=\sigma\tau'^{-1}\tau^{-1}\sigma=\sigma(\tau\tau')^{-1}\sigma=\tau\tau'$.
Hence $H$ is abelian.
For any $\rho\in G$ and $\tau\in H$, we conclude that $\rho\tau\rho^{-1}=\tau$ if $\rho\in H$ and $\rho\tau\rho^{-1}=\tau^{-1}$ otherwise.
Therefore every subgroup of $H$ is normal in $G$.
\end{proof}

\begin{lem}\label{L:extensions}
Assume that $s(K)>2$. Then the following hold:
\begin{enumerate}[$(a)$]
    \item $K(\sqrt{-1})$ is not contained in any cyclic degree-$4$ extension of $K$.
    \item If $p(K)>1$, then $K(\sqrt{-1})$ admits a dihedral extension of degree $8$.
\end{enumerate}
\end{lem}
\begin{proof}
$(a)$ This follows from \cite[Theorem 2]{Kim90}.

$(b)$
Set $M=K(\sqrt{-1})$.
We assume that $p(K)>1$ and fix  $a\in\mg{K}\setminus\sq{K}$ which is a sum of two squares in $K$.
Since $s(K)>2$, we obtain that $a\notin\sq{M}$, whereby $[M(\sqrt{a}):M]=2$.
We fix $x,y\in K$ such that $a=x^2+y^2$.
Note that $a$ is represented by the quadratic form $X^2-aY^2$ over $K$.
Since the norm of $x+y\sqrt{-1}\in M$ for the extension $M/K$ is equal to $a$, it follows by \cite[Chap.~VII, Theorem 5.10]{Lam05} that there exists some $z\in \mg{K}$ such that
$z(x+y\sqrt{-1})$ is represented by $X^2-aY^2$ over $M$.
We set $b=z(x+y\sqrt{-1})\in M$. 
Then the norm of $b$ for the extension $M/K$ is equal to $z^2a$, which is not a square in $K$.
We conclude that 
$[M(\sqrt{a},\sqrt{b}):M]=4$.
Since $b$ is represented by $X^2-aY^2$ over $M$, \cite[Theorem 5]{Kim90} yields that  $M(\sqrt{a},\sqrt{b})/M$ extends to a dihedral degree-$8$ extension of $M$.
\end{proof}

The following theorem and its proof are a variation and partial extension of \cite[Chap.~III, Theorem 1]{Bec78}. More precisely, the result there gives the equivalences between $(i)$, $(ii)$, $(v)$ and $(vi)$ of \Cref{T:Becker-hp-equiv}.
For the sequel, we essentially need the (new) implication $(i\Rightarrow iv)$, in particular to reach \Cref{C:hp-general-transc-p-bound} via \Cref{P:hp-scg}.
Our proof of $(i\Rightarrow iv)$, however, involves several of the other implications.

\begin{thm}[Becker]\label{T:Becker-hp-equiv}
Assume that $K$ is real. The following are equivalent:
\begin{enumerate}[$(i)$]
    \item $K$ is hereditarily pythagorean.
    \item Every finite nonreal  extension of $K$ contains $\sqrt{-1}$.
    \item Every finite nonreal  extension of $K$ is a Galois extension.
    \item Every finite nonreal extension of $K$ is given by $L(\sqrt{-1})$ for some finite real extension $L/K$.
    \item Every $K$-automorphism  of an algebraic closure of $K$ either fixes $\sqrt{-1}$ or has order $2$.
    \item The absolute Galois group of $K(\sqrt{-1})$ is abelian.
\end{enumerate}
\end{thm}
\begin{proof}
$(i\Rightarrow ii)$
This follows from \Cref{T:Pfister-BVG} for $n=0$. 

$(ii\Rightarrow v)$ 
We fix an algebraic closure $\Omega$ of $K$ and denote the Galois groups of $\Omega/K$ and $\Omega/K(\sqrt{-1})$ by $G$ and $H$, respectively.
Consider an arbitrary element $\sigma\in G\setminus H$.
Let $L$ denote the fixed field of $\sigma$ in $\Omega$.
Then $K\subseteq L$ and $\sqrt{-1}\notin L$. 
It follows from $(ii)$ that $L$ is real.
Let $R$ denote a maximal real subfield of $\Omega$ containing $L$.
Then $R$ is a real closure of $L$ with respect to some ordering.
In particular, $\Omega=R(\sqrt{-1})$.
By the definition of $L$, any finite extension $N$ of $L$ in $\Omega$ is cyclic with Galois group generated by $\sigma|_N$. This implies that $\sigma(R)=R$.
Since $R$ is a real closure of $L$, it follows by \cite[Theorem 1.3.2]{BCR98} that $\sigma|_R$ is the identity on $R$.
Hence $L=R$.
Therefore, the order of $\sigma$ is equal to $[\Omega:R]=2$.

$(v\Rightarrow vi)$ This is immediate by \Cref{L:groups}.

$(vi\Rightarrow i)$
Suppose that $K$ is not hereditarily pythagorean. 
Hence there exists a real finite field extension $L/K$ such that $p(L)>1$. 
It follows by \Cref{L:extensions} $(b)$ that $L(\sqrt{-1})$ admits a nonabelian Galois extension.
Therefore, the absolute Galois group of $K(\sqrt{-1})$ is not abelian.

$(ii,v\Rightarrow iii)$
This likewise follows by \Cref{L:groups}.

$(iii\Rightarrow i)$ 
Let $L/K$ be a finite real extension. 
Set $M=L(\sqrt{-1})$.
Assuming that $(iii)$ holds, every quadratic extension of $L(\sqrt{-1})$ is a Galois extension of $L$.
By \Cref{L:extensions} $(a)$, $M/L$ cannot be contained in a cyclic degree-$4$ extension of $L$.
Combining these, we conclude that every quadratic extension of $M$ is of the form $M(\sqrt{a})$ for some $a\in\mg{L}$.
In other words, we have $\mg{M}=\mg{L}\sq{M}$.
Consequently, every norm of the extension $M/L$ is a square in $L$. 
Therefore $p(L)=1$.

$(i,ii,iii\Rightarrow iv)$
Let $M/K$ be a nonreal finite extension.
By $(iii)$, $M/K$ is a Galois extension.
We choose a subextension $L/K$ of $M/K$ which is maximal for having $-1\notin \sq{L}$.
Then clearly $L$ has no proper odd-degree extension inside $M$, and $L(\sqrt{-1})$ is the unique quadratic extension of $L$ contained in $M$.
Since $M/K$ is a Galois extension, 
it follows that $M/L$ is a $2$-extension.
By $(ii)$, the field $L$ is real. By $(i)$, we therefore have $p(L)=1$.
Hence, the norm in $L(\sqrt{-1})/L$ of any element of $L(\sqrt{-1})$ is a square in $L$.
By \cite[Chap.~VII, Theorem 3.8]{Lam05}, this implies that $\mg{L(\sqrt{-1})}=\mg{L}\sq{L(\sqrt{-1})}$.
Since $L(\sqrt{-1})$ is the only quadratic extension of $L$ inside $M$, we obtain that $L(\sqrt{-1})$ has no quadratic extension inside $M$.
Since $M/L$ is a $2$-extension, we conclude that $M=L(\sqrt{-1})$.

$(iv\Rightarrow ii)$ This implication is trivial.
\end{proof}

\begin{prop}[Becker]\label{P:hp-scg}
Let $K$ be hereditarily pythagorean and let $n\in\nat^+$ be such that $|\scg{K}|=2^{n}$.
Then 
$|\scg{L}|=2^n$ for every finite real field extension $L/K$, and $|\scg{M}|=2^{n-1}$ for every finite non-real field extension $M/K$.
\end{prop}

\begin{proof}
Consider a finite real extension $L/K$. 
Then \cite[Chap.~III, Theorem 18]{Bec78} yields that $|\scg{L}|=2^n$, and for $M=L(\sqrt{-1})$, we conclude by \cite[Chap.~III, Theorem 12]{Bec78} that $|\scg{M}|=\frac{1}2 |\scg{L}|=2^{n-1}$.
In view of \Cref{T:Becker-hp-equiv}, this shows the statement.
 \end{proof}

\begin{prop}\label{P:hens-scg-hp}
Let $(K,v)$ be a henselian valued field with $\car (Kv)\neq 2$.
Then $|\scg{K}|=[vK:2vK]\cdot |\scg{(Kv)}|$.
Furthermore, $K$ is hereditarily pythagorean if and only if $Kv$ is hereditarily pythagorean.
\end{prop}

    \begin{proof}
On the one hand, the hypotheses on $(K,v)$ imply that $1+\mfm_v\subseteq\sq{\mc{O}_v}$, from which it follows that $\scg{\mc{O}_v}\simeq \scg{(Kv)}$.
On the other hand, the inclusion $\mg{\mc{O}}_v\hookrightarrow \mg{K}$ and the homomorphism $v: \mg{K} \to vK$ induce a short exact sequence $1\to \scg{\mc{O}_v}\to \scg{K}\to vK/2vK\to 1$.
These ingredients together yield that $|\scg{K}|=[vK:2vK]\cdot |\scg{(Kv)}|$.

Recall that for any finite field extension $L/K$ there is a unique extension of $v$ to a valuation $w$ on $L$.
We have that $w$ is henselian and $Lw/Kv$ is also a finite extension.
Conversely, for any finite field extension $k/Kv$, there exists a finite field extension $L/K$ such that the unique extension $w$ of $v$ to $L$ satisfies $Lw = k$.
For any such finite field extension $L/K$ with $w$ extending $v$, using that $w$ is henselian, we obtain by \Cref{inequality p'} that $L$ is real pythagorean if and only if $Lw$ is real pythagorean.
This establishes that $K$ is  hereditarily pythagorean if and only if $Kv$ is hereditarily pythagorean.
\end{proof}
    
This proposition suggests that hereditarily pythagorean fields proliferate easily via standard valuation constructions, as illustrated by the following examples.

\begin{exs}\label{EX:hp}
Let $k$ be a field. 
The field $K=k(\!(t)\!)$, called the \emph{field of Laurent series in $t$ over $k$}, carries a canonical valuation $v$, called the \emph{$t$-adic valuation}, and we have $vK=\zz$.
Moreover $(K,v)$ is complete and henselian.
The field of \emph{Puiseux series in $t$ over $k$} is given by $K'=K(\sqrt[n]{t}\mid n\in\nat^+)$. 
This is an algebraic extension of $K$, and the valuation $v$ extends uniquely to a valuation $v'$ on $K'$.  
We have that $(K',v')$ is henselian and $v'K'=\qq$.

Assume now that $k$ is hereditarily pythagorean.
Then we readily obtain by applying \Cref{P:hens-scg-hp} that $K$ and $K'$ are hereditarily pythagorean fields, too, and that $|\scg{K}|=2\cdot|\scg{k}|$ and $|\scg{K'}|=|\scg{k}|$.
\end{exs}

The torsion in powers of the fundamental ideal of the Witt ring for field extensions of a hereditarily pythagorean field $K$ depends crucially on the size of $\scg{K}$.
The following theorem partially extends \cite[Theorem]{EW87}, which covers the case $n=1$.

\begin{thm}\label{C:hp-general-transc-p-bound}
Assume that $K$ is hereditarily pythagorean and $n\in\nat^+$ is such that $|\scg{K}|=2^{n}$.
Let $m\in\nat$ and let $F/K$ be a field extension of transcendence degree $m$.
Then $\I^{n+m}F$ is torsion-free and $p(F)\leq 2^{n+m-1}$.
\end{thm}
\begin{proof}
Consider a finite field extension $M/K(\sqrt{-1})$.
By \Cref{P:hp-scg}, we have $|\scg{M}|=2^{n-1}$.
It follows by \cite[Chap.~XI, Theorem 6.4]{Lam05} that every quadratic form over $M$ of dimension larger than $2^{n-1}$ is isotropic. In particular $\I^nM=0$.
Since this holds for every finite field extension $M/K(\sqrt{-1})$, 
it follows using the exact sequence of \cite[Cor.~21.3]{EKM08} that $\I^{n+1} K(\sqrt{-1})(X)=0$.
The statement then follows by \Cref{C:transcendental-torsfree-p-bound}.
\end{proof}

Given a hereditarily pythagorean field $K$, the impact of \Cref{C:hp-general-transc-p-bound} is stronger when the group $\scg{K}$ is small. We will therefore now turn our attention to the cases where $|\scg{K}|$ is $2$ or $4$. The choice to single out precisely these two cases will become clear with \Cref{T:HP-val} below.

The field $K$ is called \emph{euclidean} if $K$ is real and $|\scg{K}|=2$, or equivalently, if $\sq{K}\cup\{0\}$ is a field ordering of $K$.
In particular, any euclidean field is pythagorean.
We call the field $K$ \emph{half-euclidean} if $K$ is real pythagorean with $|\scg{K}|=4$.
This is equivalent to $K$ being pythagorean and having precisely two field orderings.

We call $K$ \emph{hereditarily euclidean} (respectively \emph{hereditarily half-euclidean}) if $K$ is real and every finite real field extension of $K$ is euclidean (respectively half-euclidean).
In view of \Cref{P:hp-scg}, $K$ is hereditarily euclidean (respectively~hereditarily half-euclidean) if and only if $K$ is hereditarily pythagorean and further euclidean (respectively~half-euclidean).
We record explicitly two immediate consequences of \Cref{C:hp-general-transc-p-bound} for these two types of  fields, noting that they are well known in the case where $K$ is hereditarily euclidean.

\begin{cor}\label{C:hehhe-transc-p-bound}
Let $m\in\nat$ and let $F/K$ be a field extension of transcendence degree $m$.
If $K$ is hereditarily euclidean, then $p(F)\leq 2^m$. If $K$ is hereditarily half-euclidean, then $p(F)\leq 2^{m+1}$.
\end{cor}
\begin{proof}
This is \Cref{C:hp-general-transc-p-bound} for $n=1$ and $n=2$, respectively.
\end{proof}

\begin{cor}
\label{C:p-tilde-eucl-half}
If $K$ is hereditarily euclidean, then $\wt{p}(K)=3$.
If $K$ is hereditarily half-euclidean, then $\wt{p}(K)$ is either $3$ or $5$.
\end{cor}
\begin{proof}
By \Cref{P:p-tilde}~$(\ref{p-tilde-s>1})$, we have $\wt{p}(K)\geq 3$ and  $\wt{p}(K)\neq 4$.

Assume first that $K$ is hereditarily euclidean. 
By \Cref{C:hehhe-transc-p-bound}, for every function field in one variable $F/K$, we have 
$p(F)\leq 2$ and thus $p'(F)\leq 3$.
This shows that $\wt{p}(K)=3$.

Assume now that $K$ is hereditarily half-euclidean. 
By \Cref{C:hehhe-transc-p-bound}, for every function field in one variable $F/K$, we have 
$p(F)\leq 4$ and thus $p'(F)\leq 5$.
This shows that 
$\wt{p}(K)\leq 5$.
\end{proof} 

Hereditarily euclidean fields are characterised by \cite[Satz 1.2]{PZ75} as the uniquely ordered fields whose real closures are direct limits of finite extensions of odd degree.
Let us look at some ways to obtain hereditarily half-euclidean fields.

\begin{exs}\label{EX:he}
$(1)$ The field $\rr(\!(t)\!)$ is hereditarily half-euclidean; see \Cref{EX:hp}.
On this field, the $t$-adic valuation is henselian and its residue field is $\rr$, which is hereditarily euclidean.

$(2)$ Let $\ovl{\qq}$ be the algebraic closure of $\qq$ in $\cc$.
Fix a square-free integer $d>1$ and $\delta\in\rr$ with $\delta^2=d$.
Let $R_1$ and $R_2$ be two real closed subfields of $\ovl{\qq}$ containing $\delta$ and such that $\delta\in\sq{R}_1$ and $-\delta\in\sq{R}_2$.
Set $K=R_1\cap R_2$. It follows by \cite[Chap.~IV, Theorem~1]{Bec78} that $K$ is hereditarily pythagorean. 
Moreover, it is easy to see that $K$ is half-euclidean.
We conclude by \Cref{P:hp-scg} that $K$ is hereditarily half-euclidean. 

In this case \Cref{T:HP-val} applies trivially, namely by taking the trivial valuation on $K$, and in fact this is the only possibility, as $K$ is an algebraic extension of $\qq$ and therefore admits no non-trivial valuation with a real residue field. 
Therefore hereditarily half-euclidean residue fields cannot be omitted in \Cref{T:HP-val}.

$(3)$ Let $k$ be any hereditarily half-euclidean field, e.g.~one of the fields in $(1)$ and $(2)$. Then the field of Puiseux series $k(\!(t)\!)(\sqrt[n]{t}\mid n\in\nat^+)$ from \Cref{EX:hp} is hereditarily half-euclidean as well.
\end{exs}

Note that only in \Cref{EX:he}~$(2)$ did we see a non-trivial example of a hereditarily pythagorean field that did not involve a henselian valuation in the construction.
The following characterisation due to L.~Br\"ocker \cite{Bro76}
of hereditarily pythagorean fields indicates that this is not at all a coincidence and underlines the crucial role of hereditarily euclidean and hereditarily half-euclidean fields in the study of general hereditarily pythagorean fields.

\begin{thm}[Br\"ocker]
\label{T:HP-val}
A field $K$ is hereditarily pythagorean if and only if there exists a henselian valuation on $K$ whose residue field is hereditarily euclidean or hereditarily half-euclidean.
\end{thm}

\begin{proof}
The statement is contained in \cite[Prop.~3.5]{Bro76}.
\end{proof}

Now we will see how helpful this theorem is for the study of Pythagoras numbers of function fields in one variable over hereditarily pythagorean fields.

\begin{thm}\label{C:final}
Assume that $K$ is hereditarily pythagorean.
Then $p(F)\leq 5$ for every function field in one variable $F/K$.
Furthermore, $p(K(X,Y))\leq 8$.
\end{thm}
\begin{proof}
By \Cref{T:HP-val}, there exists a henselian valuation $v$ on $K$ such that $Kv$ is hereditarily euclidean or hereditarily half-euclidean.
In particular, $Kv$ is real, whereby $\car(Kv)=0$.
We conclude by \Cref{T:p-tilde-henselian-equality} and \Cref{C:p-tilde-eucl-half}
that $\wt{p}(K)=\wt{p}(Kv)\leq 5$, and then by \Cref{P:p-tilde}~$(\ref{p-tilde-rafufi2})$ that $p(K(X,Y))\leq 8$.
\end{proof}

\begin{rem}\label{R:final}
If $-1\in\sq{K}$ then $p(F)\leq 2$ for any field extension $F/K$. Hence \Cref{T:main-2}
is covered by \Cref{C:final}, in view of \Cref{T:Becker}.
\end{rem}

\begin{qu}\label{Q:final1}
Which implications hold between the following conditions?
\begin{enumerate}[$(i)$]
    \item $K$ is hereditarily pythagorean.
    \item $-1\notin\sq{K}$ and $p(F)=2$ for some regular function field $F/K$.
    \item $K$ is real and $\wt{p}(K)=3$.
    \item $p(K(X,Y))=4$.
\end{enumerate}
\end{qu}
Certainly $(i\Rightarrow ii)$ holds, as one may take $F=K(X)$, by \Cref{T:Becker}.
Whether the converse implication holds was asked in \cite[Question~4.4]{BVG09}.
\Cref{C:fufi2-p4} establishes the equivalence $(iii\Leftrightarrow iv)$.
By means of \Cref{T:HP-val} and \Cref{T:p-tilde-henselian-equality}, the problem of the validity of $(i\Rightarrow iii, iv)$ can be reduced to the case where $K$ is either hereditarily euclidean or hereditarily half-euclidean. 
When $K$ is hereditarily euclidean, the positive answer follows by \Cref{C:p-tilde-eucl-half}.
So the crucial question remaining in this context is the following.

\begin{qu}\label{Q:final2}
Is $\wt{p}(K)=3$ whenever $K$ is hereditarily half-euclidean?
\end{qu}

\subsection*{Acknowledgments} 
We want to thank Eberhard Becker and David Leep for inspiring discussions and very helpful comments on the topic and its history, and Eva Bayer-Fluckiger for stimulating us to improve the presentation.
We are grateful to the anonymous referee for a very careful reading and various valuable comments helping us to improve the clarity of the exposition.

We gratefully acknowledge support from the {FWO Odysseus programme} (pro\-ject G0E6114N, {Explicit Methods in Quadratic Form Theory}), 
the Bijzonder Onderzoeksfonds (BOF), Universiteit Antwerpen
(project BOF-DOCPRO4, 2533), CONICYT-PFCHA/Doctorado Nacional/2017-folio 21170477, the Universidad de Santiago de Chile (Proyecto DICYT,
c\'odigo 042432G), ANID Fondecyt\linebreak Postdoctoral Grant 3240191, ANID Concurso de Subvenci\'on a la Instalaci\'on en la Academia 2025 N$^\circ$85250089, the {Czech Science Foundation} (GA\v CR) grant 21-00420M, and finally the {Charles University Research Centre programmes} \linebreak UNCE/SCI/022 and PRIMUS/24/SCI/010.

\subsection*{Conflict of interest statement}
On behalf of all authors, the corresponding author states that there is no conflict of interest.

\subsection*{Data accessibility statement}
The  article describes entirely theoretical research.
All new data supporting the findings presented in this article are included in the manuscript.
Data sharing is not applicable to this article.

\subsection*{AI usage statement}
The authors declare that they used AI tools (\emph{ChatGPT}, \emph{Google Gemini}, \emph{Claude}) during the finalization process for this article to identify inconsistencies and lack of clarity.

\bibliographystyle{plain}

\end{document}